\documentclass[10pt,a4paper]{article}
\usepackage[latin1]{inputenc}
\usepackage{amsmath}
\usepackage{amsfonts}
\usepackage{amssymb}
\usepackage{wasysym}
\usepackage{graphicx}
\usepackage{geometry}
\usepackage{amsmath}
\usepackage{amssymb}
\usepackage{amsthm}
\usepackage{color}
\usepackage{tikz}
\newtheorem{mydef}{Definition}
\newtheorem{mylem}{Lemma}
\newtheorem{mythm}{Theorem}
\DeclareMathOperator{\Hom}{Hom}
\DeclareMathOperator{\Inj}{Inj}
\title{Rational exponents for hypergraph Turan problems}
\author{Matthew Fitch \footnote{Zeeman Building,
University of Warwick
Coventry CV4 7AL, United Kingdom. Research supported by ERC Grant No. 306493. Email address: M.H.D.Fitch@warwick.ac.uk}}
\begin{document}
\maketitle

\paragraph*{Abstract\\}

Given a family of $k$-hypergraphs $\mathcal{F}$, $ex(n,\mathcal{F})$ is the maximum number of edges a $k$-hypergraph can have, knowing that said hypergraph has $n$ vertices but contains no copy of any hypergraph from $\mathcal{F}$ as a subgraph. We prove that for a rational $r$, there exists some finite family $\mathcal{F}$ of $k$-hypergraphs for which $ex(n,\mathcal{F})=\Theta(n^{k-r})$ if and only if $0\leq r \leq k-1$ or $r=k$.\\

\section{Introduction}

\paragraph{Definitions:}

Given an integer $k\geq 2$, a \emph{$k$-hypergraph} $G$ is a set of points (called the \emph{vertices}), together with a collection of $k$-subsets of the vertices (called the \emph{edges}). For such a $k$-hypergraph, $|G|$ refers to its number of vertices and $e(G)$ refers to the number of edges.\\

Given $k$-hypergraphs $G$ and $X$, a \emph{graph homomorphism} (often shortened to \emph{homomorphism}) from $G$ to $X$ means a function $f$ that assigns to each vertex of $G$ some vertex in $X$ and that also preserves edges, ie for every edge $\{x_1,x_2...,x_k\}$ in $G$, $\{f(x_1),f(x_2),...,f(x_k)\}$ is also an edge of $X$. The set of all such graph homomorphisms is called $\Hom(G,X)$. The image of such a graph homomorphism is called a \emph{homomorphic copy} of $G$ in $X$.\\

Given $k$-hypergraphs $G$ and $X$, $\Inj(G,X)$ is a subset of $\Hom(G,X)$, and is defined to be the set of all injective graph homomorphisms from $G$ to $X$, ie, those homomorphisms $f$ with the property that for any 2 distinct vertices $x$ and $y$, $f(x) \neq f(y)$. When $|\Inj(G,X)|\geq 1$ we say that $X$ \emph{contains $G$ as a subgraph}.\\

Given an integer $k\geq 2$ and a family $\mathcal{F}$ of $k$-hypergraphs, $ex(n,\mathcal{F})$ is defined to be the maximum number of edges across all $k$-hypergraphs that have $n$ vertices and do not contain any element of $\mathcal{F}$ as a subgraph. In general, this quantity can be anything from $0$ (as in the case $\mathcal{F}=\{E\}$, where $E$ is just a single edge) to $\binom{n}{k}$ (as in the case where $\mathcal{F}$ is empty). \\

Finding $ex(n,\{F\})$ for a fixed graph or hypergraph $F$ is known as the Tur\' an problem. For ordinary ($k=2$) non-bipartite graphs, we have a reasonable understanding: Tur\' an gave an exact solution when $F$ is a complete graph \cite{Turan}, while Erd\" os and Stone gave an asymptotic solution for any non-bipartite graph \cite{ErdosStone}. However, for bipartite graphs and more general hypergraphs ($k\geq 3$), very little is known, not even asymptotically. It is a major open problem in extremal combinatorics to come up with some sort of understanding of these numbers. \\

For a lot of families of hypergraphs, $ex(n,\mathcal{F})$ is of order $\Omega(n^k)$. However, there are some for which $ex(n,\mathcal{F})$ is of order $o(n^k)$. We call this case a \emph{degenerate Tur\' an problem}. This categery includes the bipartite graph case and these are going to be the ones we will look at in this paper.\\

In 1979, Erd{\H{o}}s conjectured that for every rational $r$ between 1 and 2, there exists a finite family of bipartite graphs $\mathcal{F}$ with $ex(n,\mathcal{F})=\Theta(n^{r})$  \cite{Erdos}. This conjecture was later proved in 2015 by Bukh and Conlon \cite{BukhConlon}.  \\

In 1986, Frankl proved a related result for hypergraphs: for every rational $r\geq 1$, there exists some $k\in\mathbb{N}$ and some finite family $\mathcal{F}$ of $k$-hypergraphs such that $ex(n,\mathcal{F})$ is of order $n^r$. (Sidenote: the $\mathcal{F}$ that Frankl used also had the property that every $F\in \mathcal{F}$ had exactly 2 edges.) \cite{Frankl}.\\

In 2016, Ma, Yuan and Zhang discovered an infinite family of $k$-hypergraphs for which they could solve the Tur\' an problem asymptotically. They proved that $K^{(k)}_{s_1,s_2,...,s_k}$, the complete $k$-partite $k$-hypergraph with partition sizes $s_1,s_2,...,s_k$ has $ex(n,K^{(k)}_{s_1,s_2,...,s_k})=\Theta(n^{k-\frac{1}{s_1s_2s_3...s_k}})$ \cite{MaYuanZhang}.\\

In this paper, we will extend Bukh and Conlon's result to hypergraphs, ie, we will prove that for every rational $r$ between 1 and $k$, there exists a finite family of $k$-hypergraphs $\mathcal{F}$ with $ex(n,\mathcal{F})=\Theta(n^{r})$. This also an improvement on Frankl's result since we now have a family of $k$-hypergraphs for all $k\geq r$, instead of for just one specific $k$. This is also of interest as an infinite family of $k$-hypergraphs for which the answer to the Tur\' an problem is known.\\

To prove this we will use similar methods as Bukh and Conlon, both in the construction of $\mathcal{F}$ and for the lower bound. However, the proof of the upper bound ($ex(n,\mathcal{F})\leq c \cdot n^{r}$ for some constant $c$) does not easily generalise to hypergraphs. We come up with an alternative proof, where ~\cite{Szegedy} ends up being very helpful.\\

For simplicity, we will be exchanging $r$ and $k-r$ for the rest of the paper. \\

We will at first only consider the case where $0\leq r < 1$:

\begin{mythm}

Given a rational $r$, $0\leq r < 1$, there exists some finite collection of $k$-hypergraphs $\mathcal{F}$ such that $ex(n,\mathcal{F})=\Theta(n^{k-r})$.

\end{mythm}

 Our first section deals with the construction of the family of hypergraphs $\mathcal{F}$ that will solve Theorem 1. They will be hypergraph versions of the graphs from ~\cite{BukhConlon}. \\
 
 In our second section, we prove the lower bound, ie that $ex(n,\mathcal{F})\geq \Theta(n^{k-r})$. This involves constructing a hypergraph with $n$ vertices and $\Theta(n^{k-\frac{a}{b}})$ edges but that does not contain any copy of any hypergraph from $\mathcal{F}$.  The proof is again adapted from ~\cite{BukhConlon}.\\

 In our third section, we prove the upper bound, ie that $ex(n,\mathcal{F})\leq \Theta(n^{k-r})$. However, unlike in the first two sections, the proof from ~\cite{BukhConlon} cannot be easily extended to hypergraphs. We instead use a partial version of the Sidorenko Conjecture ~\cite{Sidorenko} (also conjectured by Simonovits ~\cite{Simonovits}) from the paper of Szegedy ~\cite{Szegedy}. When $n$ is a sufficiently large interger, this allows us to find some copy of an element of $\mathcal{F}$ in any hypergraph $X$ with $n$ vertices and with at least $\Theta(n^{k-r})$ edges, thereby proving the upper bound.\\
 
In our final section, we consider what happens for other $r$s. We first extend the result from $0\leq r < 1$ to $0 \leq r \leq k-1$:
 
 \begin{mythm}

Given a rational $r$, $0\leq r < k-1$, there exists some collection of $k$-hypergraphs $\mathcal{F}$ such that $ex(n,\mathcal{F})=\Theta(n^{k-r})$.

\end{mythm}

Observation: The case where $k-1<r<k$ is impossible. This is a corrolery of the Sunflower Lemma \cite{erdos1960intersection}, which involves hypergraphs called \emph{sunflower}s. A sunflower is a $k$-hypergraph which contains a specific \emph{kernal}, a set between 0 and $k-1$ points, and any two edges of the sunflower intersect in exactly the kernal. The sunflower lemma states that whenever $F$ is a collection of $k$-hypergraphs such that for all $0\leq i\leq k-1$, $F$ contains a sunflower with kernal size $i$, then $ex(n,F)$ is bounded by a constant (independent of $n$). We shall provide more details in the final section.\\

\paragraph{Notation:\\}

$\, \, \, \mathbb{P}(A)$ means the probability that event $A$ will occur.\\

$\mathbb{E}(B)$ means the expected value of the random variable $B$\\

\subsection*{Algebraic Geometry}

The proof will use some algebraic geometry. What follows in this section is a brief overview of the results we will use. See ~\cite{alggeom} for more information and proofs.

\begin{mydef} (1.1.2 in \cite{alggeom}) Given an algebraically closed field $\mathbb{F}$, an \emph{affine algebraic variety} $V$ (often shorted to just variety) over $\mathbb{F}$ is a set of the form:\\ 
$V= \{(x_1,x_2,...,x_n) \in \overline{\mathbb{F}}^n \, | \, P_1(x_1,x_2,...,x_n)=P_2(x_1,...,x_n)= ... = P_m(x_1,...,x_n)=0\}$, where $P_1,P_2,...,P_m$ are polynomials over $\mathbb{F}$ with $n$ variables. \\

\end{mydef}

\begin{mylem} (1.1.4 and 1.1.5 in \cite{alggeom}) If $U$ and $V$ are varieties over $\mathbb{F}$,  then $U\cap V$ and $U\cup V$ are also varieties over $\mathbb{F}$.\\

\end{mylem}

\begin{mydef} (1.1.10 in \cite{alggeom}) Given a variety $V$ over $\mathbb{F}$, we say that $V$ is \emph{reducible} if there exist varieties $U,U'\subsetneq V$ such that $V=U\cup U'$. If $V$ is not reducible, we say it is \emph{irreducible}.

\end{mydef}

\begin{mylem} (1.1.12a in \cite{alggeom}) A variety $V$ can be decomposed uniquely (up to ordering) into maximal irreducible components: $V=U_1 \cup U_2 \cup ... \cup U_k$ where the $U_i$ are all irreducible. \\
That means that if we have two decompositions $V = \bigcup_{i=1}^k U_i = \bigcup_{j=1}^l U'_j$, then for every $1\leq i \leq k$, there exists some $1\leq j \leq l$ such that $U_i \subseteq U'_j$ and vice-versa. \\
\\
Furthermore, the number of components in such a decomposition is bounded above by $d^m$ where $m$ is the number of polynomials that generate the variety, and $d$ is their maximum degree.

\end{mylem}

\begin{mydef} (1.2.15 to 1.2.17 in \cite{alggeom}) Given a non-empty irreducible variety $V$ over $\mathbb{F}$, its \emph{dimension} $\delta$ is the length of the longest sequence: $V=V_\delta\supsetneq V_{\delta-1} \supsetneq V_{d_2} \supsetneq ... \supsetneq V_0 \supsetneq \emptyset$, where every $V_i$ is irreducible. This is well defined for every non-empty irreducible variety.\\
When $V$ is reducible, we say its dimension is the largest dimension of one of its irreducible components.

\end{mydef}

It is fairly easy to see that a finite set of points has dimension 0, the space $\overline{\mathbb{F}}^n$ has dimension $n$, and that when $V$ is a non-empty variety generated by $k$ polynomials in $n$ variables: $P_1,P_2,...,P_m \in \mathbb{F}[X_1,X_2,...,X_n]$, then $V$ has dimension at least $n-m$.\\
\\

Although we require $\mathbb{F}$ to be algebraically closed for the theory to work, most practical applications involve fields that are not algebraically closed. However, this isn't a problem because if $\mathbb{F}'$ is an arbitrary field, then it has an algebraic closure $\overline{\mathbb{F}'}$. We can then the use properties of algebraic varieties over $\overline{\mathbb{F}'}$ to say things about the corresponding set over $\mathbb{F}'$:

\begin{mydef} Given a variety $V$ over an algebraically closed field $\mathbb{F}$ and a subfield $\mathbb{F}'\subseteq \mathbb{F}$ (which might not be algebraically closed), the $\mathbb{F}'$-\emph{rational points} of the variety, denoted by $V(\mathbb{F}')$, are defined to be the points of $V$ that can be written using elements of $\mathbb{F}'$, ie: $V(\mathbb{F}') = V\cap \mathbb{F}'^n$.

\end{mydef}

\begin{mythm} (Lang-Weil bound) ~\cite{LangWeil} Let $\mathbb{F}_p$ be the finite field of order $p$, where $p$ is a power of a prime. Let $V$ be an irreducible variety of dimension $\delta$ over $\overline{\mathbb{F}_p}$. Then $V(\mathbb{F}_p)$ is either empty or has $|V(\mathbb{F}_p)| = p^\delta(1+O(p^{-1/2}))$ 

\end{mythm}

\section{The set of hypergraphs}

Since $r$ is a rational smaller than 1, we can write $r=\frac{a}{b}$ for $a$ and $b$ positive integers and $b>a$. By multiplying $a$ and $b$ by some constant, we can assume without loss of generality that $b\geq a-k+3$. Now given $a$,$b$,$k$ integers such that $b\geq a-k+3$, consider the hypergraph as in the picture:


\begin{tikzpicture}[scale=0.17]

\draw (0,0) node {$\circ$} ;
\draw (10,0) node {$\circ$} ;
\draw (20,0) node {$\circ$} ;
\draw (30,0) node {$\circ$} ;
\draw (40,0) node {$\circ$} ;
\draw (50,0) node {$\circ$} ;
\draw (56,0) node {...} ;
\draw (62,0) node {$\circ$} ;

\draw (5,5) node {$\bullet$};
\draw (5,10) node {$\bullet$};
\draw (5,15) node {$\bullet$};
\draw (15,5) node {$\bullet$};
\draw (15,10) node {$\bullet$};
\draw (15,15) node {$\bullet$};
\draw (25,5) node {$\bullet$};
\draw (25,10) node {$\bullet$};
\draw (35,5) node {$\bullet$};
\draw (35,10) node {$\bullet$};
\draw (35,15) node {$\bullet$};
\draw (45,10) node {...};

\draw[color=green] (-1.5,-0.5) -- (11.5,-0.5);
\draw[color=green] (-1.5,-0.5) -- (5,6);
\draw[color=green] (11.5,-0.5) -- (5,6);
\draw[color=green] (-2,-1) -- (12,-1);
\draw[color=green] (-2,-1) -- (5,11.5);
\draw[color=green] (12,-1) -- (5,11.5);
\draw[color=green] (-2.5,-1.5) -- (12.5,-1.5);
\draw[color=green] (-2.5,-1.5) -- (5,17);
\draw[color=green] (12.5,-1.5) -- (5,17);

\draw[color=green] (8.5,-0.75) -- (21.5,-0.75);
\draw[color=green] (8.5,-0.75) -- (15,6);
\draw[color=green] (21.5,-0.75) -- (15,6);
\draw[color=green] (8,-1.25) -- (22,-1.25);
\draw[color=green] (8,-1.25) -- (15,11.5);
\draw[color=green] (22,-1.25) -- (15,11.5);
\draw[color=green] (7.5,-1.75) -- (22.5,-1.75);
\draw[color=green] (7.5,-1.75) -- (15,17);
\draw[color=green] (22.5,-1.75) -- (15,17);

\draw[color=green] (18.5,-0.5) -- (31.5,-0.5);
\draw[color=green] (18.5,-0.5) -- (25,6);
\draw[color=green] (31.5,-0.5) -- (25,6);
\draw[color=green] (18,-1) -- (32,-1);
\draw[color=green] (18,-1) -- (25,11.5);
\draw[color=green] (32,-1) -- (25,11.5);

\draw[color=green] (28.5,-0.75) -- (41.5,-0.75);
\draw[color=green] (28.5,-0.75) -- (35,6);
\draw[color=green] (41.5,-0.75) -- (35,6);
\draw[color=green] (28,-1.25) -- (42,-1.25);
\draw[color=green] (28,-1.25) -- (35,11.5);
\draw[color=green] (42,-1.25) -- (35,11.5);
\draw[color=green] (27.5,-1.75) -- (42.5,-1.75);
\draw[color=green] (27.5,-1.75) -- (35,17);
\draw[color=green] (42.5,-1.75) -- (35,17);

\draw[color=red] (10,0) ellipse (12 and 2);
\draw[color=red] (20,0) ellipse (12 and 2);
\draw[color=red] (30,0) ellipse (12 and 2);
\draw[color=red] (40,0) ellipse (12 and 2);

\draw (65,10) node {\Huge{\}} \LARGE{$b-a+k-1$}};
\draw (65,0) node {\LARGE{\} $a$}};

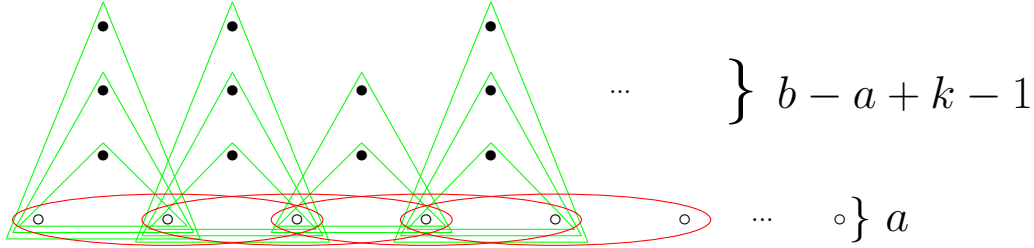
\draw (40,-5) node {Example of the hypergraph $\mathcal {T}$ in the case $k=3$};

\end{tikzpicture}

It is essentially a hypergraph version of the graph from ~\cite{BukhConlon}. It is comprised of an ordered set of $a$ vertices (in white) with edges (in red) being sets of $k$ vertices in a row. We add to this $b-a+k-1$ vertices (in black) and for each one, an edge (in green) connecting it to $k-1$ vertices in a row. This makes the total number of edges to be $b$. These black vertices are as evenly spaced as possible (see picture below). Formally, the $i$th black vertex is connected to the $\lfloor 1+\frac{(i-1)(a-k+2)}{b-a+k-2}\rfloor$th $(k-1)$-set of consecutive white vertices. There is one exception, and that is the last (ie: $(b-a+k-1)$th) black vertex is connected to the last (ie: $(a-k+2)$th) consecutive set of white vertices, not, as the formula suggests, the $(a-k+3)$th, because that one doesn't exist. We call the vertex-set of this hypergraph $T$, the subset of black vertices $R$ and call these black vertices the \emph{roots} of $\mathcal{T}$.\\
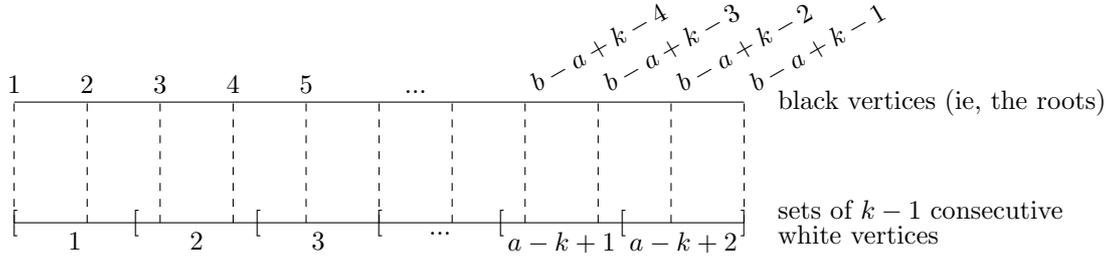
\noindent
\begin{tikzpicture}[scale=0.8]

\draw (0,0) -- (12,0);
\draw (0,2) -- (12,2);

\draw (12.4,0.2) node[right] {sets of $k-1$ consecutive};
\draw (12.4,-0.2) node[right] {white vertices};  
\draw (12.4,2) node[right] {black vertices (ie, the roots)}; 

\draw (0,0) node {[};
\draw (2,0) node {[};
\draw (4,0) node {[};
\draw (6,0) node {[};
\draw (8,0) node {[};
\draw (10,0) node {[};
\draw (12,0) node {]};

\draw (1,-0) node[below] {$1$};
\draw (3,-0) node[below] {$2$};
\draw (5,-0) node[below] {$3$};
\draw (7,-0) node[below] {$...$};
\draw (9,-0) node[below] {$a-k+1$};
\draw (11,-0) node[below] {$a-k+2$};

\draw[dashed] (0,0) -- (0,2);
\draw[dashed] (1.2,0) -- (1.2,2);
\draw[dashed] (2.4,0) -- (2.4,2);
\draw[dashed] (3.6,0) -- (3.6,2);
\draw[dashed] (4.8,0) -- (4.8,2);
\draw[dashed] (6,0) -- (6,2);
\draw[dashed] (7.2,0) -- (7.2,2);
\draw[dashed] (8.4,0) -- (8.4,2);
\draw[dashed] (9.6,0) -- (9.6,2);
\draw[dashed] (10.8,0) -- (10.8,2);
\draw[dashed] (12,0) -- (12,2);

\draw (0,2) node[above] {$1$};
\draw (1.2,2) node[above] {$2$};
\draw (2.4,2) node[above] {$3$};
\draw (3.6,2) node[above] {$4$};
\draw (4.8,2) node[above] {$5$};
\draw (6.6,2) node[above] {$...$};
\draw (8.4,2.2) node[right,rotate=30] {$b-a+k-4$};
\draw (9.6,2.2) node[right,rotate=30] {$b-a+k-3$};
\draw (10.8,2.2) node[right,rotate=30] {$b-a+k-2$};
\draw (12,2.2) node[right,rotate=30] {$b-a+k-1$};

\draw (-0.3,-1.5) node[right] {\underline{An example of how the roots are connected  to the sets of $k-1$ consecutive non-roots}};
\draw (-0.3,-2.3) node[right] {In this picture for example, the second $(k-1)$-set of non-roots is connected to both the 3rd};
\draw (-0.3,-2.7) node[right] {and 4th root but no others. When a root lands exactly on a border,  it gets connected to the};
\draw (-0.3,-3.1) node[right] {$(k-1)$-set corresponding to the interval on its right EXCEPT for the very last one, which};
\draw (-0.3,-3.5) node[right] {gets connected to the $(k-1)$-set on its left (because there is nothing to the right)};

\end{tikzpicture}

\subsection{$\mathcal{T}$ is balanced\\}

\begin{mydef}

Given a set $S$ of non-roots, define $\epsilon(S)$ to be the number of edges that contain a point of $S$. 

\end{mydef}

\begin{mydef}
A rooted $k$-hypergraph $\mathcal{U}$ with vertex set $U$ and set of roots $R$ is \emph{balanced} if for any subset $S\subset U-R$, we have:\\ $\frac{\epsilon(S)}{|S|}\geq \frac{\epsilon(U-R)}{|U-R|}$. 
\end{mydef}

Notice that in the case of our hypergraph $\mathcal{T}$ (whose vertex set is $T$), we have $\epsilon(T-R)$ is the total number of edges in the hypergraph, ie $\epsilon(T-R)=b$.

\begin{mylem}

The hypergraph $\mathcal{T}$ defined above is balanced.

\end{mylem}

\paragraph*{Proof:} First of all, if every edge of $\mathcal{T}$ contains an element of $S$, then the result is trivial since $|S|\leq |T-R|$. Without loss of generality, we can therefore assume that there is at least one edge that doesn't contain any element of $S$, which means there is a section of length at least $k-1$ that does not contain any element of $S$. Call this a \emph{hole}. We also seperate $S$ into a sequence of \emph{blocks}, by which we mean a maximal sequence of elements of $S$ with no gaps between them.\\

Suppose we have a block directly to the left of a hole and also suppose that it does not contain the leftmost vertex of $T$. Call this block $R$. What happens if we shift the entire block to the left? Because the black vertices (roots) are evenly distributed, the number of green edges (edges containing roots) adjacent to $R$ varies by at most 1. The number of red edges (edges not containing roots) containing a point of $R$ stays the same unless we are reaching the left edge of $\mathcal{T}$, in which case it goes down. If we do not reach the left side of $T$, then that means there is another block in the way. In this case, the edges containing points of that block and the edges containing points of $R$ will start to coincide. Regardless of which case we are in, when we do this step, the number of red edges containing elements of $S$ goes down by at least 1, while the number of green edges changes by at most 1; therefore, the overall number of edges containing elements of $S$ goes down or stays constant, while $|S|$ stays constant. Therefore we can assume without loss of generality that this step has been completed.\\

By repeating this step multiple times, we can move blocks left until they merge with other blocks, and then continue moving the bigger blocks until we eventually have everything to the left of the hole is in one big block as left as it can go. By a similar argument, everything to the right of the hole is in one big block as far right as it goes. Say the big block on the left has size $x$ and the one on the left has size $y$.\\

If the total number of vertices in left big block is $x$ then we get $x$ red edges. The green edges we get are those that connect to the first $x$ $(k-1)$-sets. Recall from the definition, that the $ith$ green edge connects to the $\lfloor 1+\frac{(i-1)(a-k+2)}{b-a+k-2}\rfloor$th $(k-1)$-set. Therefore the number of green edges is the maximal $i$ such that $\lfloor 1+\frac{(i-1)(a-k+2)}{b-a+k-2}\rfloor \leq x$,  ie s.t. $\frac{(i-1)(a-k+2)}{b-a+k-2} < x$, ie $i=\lceil\frac{x(b-a+k-2)}{a-k+2}\rceil$. Similarly, we can calculate the number of red edges in the right big block as $y$ and the number of green edges in it as $\lfloor\frac{y(b-a+k-2)}{a-k+2}\rfloor+1$. Therefore $\epsilon(S)$ is at least $|S|+\lceil\frac{|S|(b-a+k-2)}{a-k+2}\rceil=\lceil|S|\frac{b}{a-k+2}\rceil \geq |S|\frac{b}{a}$. Thus, $\mathcal{T}$ is indeed a balanced rooted hypergraph.

\begin{flushright}
$\Box$
\end{flushright}

\begin{mydef}

$\mathcal{T}^{\leq s}$, the \emph{$s$th power of the rooted hypergraph $\mathcal{T}$}, is defined to be the set of all $k$-hypergraphs that are formed by taking the union of $s$ copies of $\mathcal{T}$ and making them agree on the roots. For the non-roots (ie, the $s$ copies of each non-root), any disposition is allowed: they can be distinct, they can coincide with each other, or they can even coincide with different non-roots from other copies of $\mathcal{T}$.

\end{mydef}

\begin{mydef}
$\mathcal{T}^{s}$=$\mathcal{T}^{\leq s} \backslash \mathcal{T}^{\leq s-1}$.
 \end{mydef}

\begin{mylem}
For any $H$ in $\mathcal{T}^s$, the number of edges in $H$ is at least $(|H|-|R|)\frac{b}{a}$. 
\end{mylem}

\paragraph*{Proof:} We prove this lemma by induction. The case $s=1$ is trivial since then $H=\mathcal{T}$.\\
 
Given $H\in \mathcal{T}^s$, we can write $v(H)$ as $v(H') \cup S$ where $H'$ is in $\mathcal{T}^{s-1}$ and $S$ is all the extra vertices from the $s$th copy of $\mathcal{T}$ that aren't already included in $H$. We can consider $S$ as a subset of $T-R$. Since $\mathcal{T}$ is balanced, we have that the number of edges containing an element of $S$ is at least $|S|\frac{b}{a}$. By induction, the number of edges in $H'$ is at least $(|H'|-|R|)\frac{b}{a}$. Therefore the total number of edges in $H$ is at least $(|S|+|H'|-|R|)\frac{b}{a}=(|H|-|R|)\frac{b}{a}$.\\

Therefore by induction, we have proved that $H$ has at least $|H-R|\frac{b}{a}$ edges.\\

\begin{flushright}
$\Box$
\end{flushright}

The set of hypergraphs $\mathcal{F}$ we will take to prove Theorem 1 is $\mathcal{T}^p$ for $a$ and $b$ such that $r=\frac{a}{b}$ and some sufficiently large $p$ (which we will define later as a function of $a$, $b$ and $k$).

\section{The lower bound}

To prove that $ex(n,\mathcal{F})\geq \Theta(n^{k-r})$, we need to construct a hypergraph $X$ with at least $\Theta(n^{k-r})$ edges but without any copies of $\mathcal{F}$. The hypergraph we will take will also be a hypergraph version of the graph from ~\cite{BukhConlon}. \\
\\
Pick constants $s= b(b-a+k-1)+a+1$ and $d= bs-1=b^2(b-a+k-1)+ab+b-1$. Then pick a sufficiently large prime $q$. In particular, we will require $q\geq \binom{d+1}{2}$.\\

The set of vertices of $X$ is $\mathbb{F}_q^b \sqcup \mathbb{F}_q^b \sqcup ... \sqcup \mathbb{F}_q^b$, where there are $k$ copies of $\mathbb{F}_q^b$. Also pick uniformly independently at random $a$ polynomials in $k$ variables of degree at most $d$: $f_1, f_2, ..., f_a$: $\mathbb{F}_q^b \times \mathbb{F}_q^b \times \mathbb{F}_q^b \times ... \times \mathbb{F}_q^b \rightarrow \mathbb{F}_q$ (there are $k$ copies of $\mathbb{F}_q^b$). [Note: picking a polynomial of degree at most $d$ at random here means that for every coefficient of degree $\leq d$, pick an element of $\mathbb{F}_q$ uniformly at random and independently of the others.] The edges of $X$ are defined to be $(x_1,x_2,...,x_k)$ such that $f_i(x_1,x_2,...,x_k)=0$ for all $i$.\\

 Thus $X$ is $k$-partite and has $kq^{b}=N$ vertices. The set of edges of $X$ are equivalent to the rational points of variety $V(\mathbb{F}_q)$, defined by $a$ polynomials: $f_1,f_2,...,f_a$. By the Lang-Weil bound ~\cite{LangWeil}, this variety is either empty or has size at least $(c-\Theta(q^{-1/2})*q^{dim(V)}$, where $c$ is some constant depending only on $k$,$a$,$b$ and $d$ but is independent of $q$. Since we have only $a$ polynomials defining the variety, we have $dim(V)\geq bk-a$ (unless its empty). Therefore, either there are 0 edges in $X$, or there are at least $\Theta(q^{bk-a})=\Theta(N^{k-\frac{a}{b}})$ edges in $X$, no matter which $f_i$s we choose.\\
 \\
\paragraph{Probability that $X$ is empty\\} 

Suppose we have already picked all the non-constant coefficients of all the $f_i$s. Pick some points $(x_1,x_2,...,x_k)$ arbitrarily. Then for each $f_i$, there is exactly one value for the constant coefficient that makes $f_i(x_1,...,x_k)=0$. The probability we pick it is $1/q$. Multiplying these together, the probability we pick exactly the right value for every $f_i$ is $1/q^a$ because we picked the functions independently of each other. Therefore $X$ contains $(x_1,x_2,...,x_n)$ (and in particular, is non-empty) with probability at least $1/q^a$. For the next parts, we'll only be considering the case where $X$ is indeed non-empty. 

\paragraph{Proof that this hypergraph is $\mathcal{T}^p$-free with positive probability\\}

Given a copy $A$ of a hypergraph $H\in\mathcal{T}^{\leq s}$ in $X$, we know it has an ordered set of $b-a+k-1$ roots. We'll call this ordered set $r(A)=(w_1,w_2,...,w_{b-a+k-1})=\textbf{w}$.\\

Before finding a suitable $A\in \mathcal{T}^m$, we'll start by picking out a potential candidate for $r(A)$. This means we are arbitrarily picking an ordered set of $b-a+k-1$ vertices: $(w_1,w_2,...,w_{b-a+k-1})=\textbf{w}$. Now in some cases, some $w_i$s might be in the wrong parts which makes it impossible for any copy of $\mathcal{T}$ to appear with those roots. We will assume we are not in this case, and that the $w_i$s are all in the correct parts so that copies of $\mathcal{T}$ are in fact possible. We will consider these $w_i$s as elements of $\mathbb{F}_q^b$.\\

Let $C$ be the set of all copies of $\mathcal{T}$ in $X$ that have $w_1,w_2,...,w_{b-a+k-1}$ as its roots. We are interested in this because whenever we have a copy of a hypergraph of $\mathcal{T}^p$ with the given roots, that implies $|C|\geq p$. For the moment, our goal will be to find an upper bound for $\mathbb{P}(|C|\geq p)$, since that will also be an upper bound on the probability of getting a copy of $\mathcal{T}^p$. 

\begin{mylem} Given $k$, $a$,$b$ and $d$, there exists some $p$ such that for all $q$ sufficiently large, $|C|\geq p \Leftrightarrow |C| \geq q/2$ 

[Note: This is how we define the $p$ when we did $\mathcal{F}=\mathcal{T}^p$. Since $d$ is defined as a function of $a$ and $b$, this $p$ only depends on $a$, $b$ and $k$.]

\end{mylem}

\paragraph*{Proof:} We will treat vertices of our hypergraph as elements in $\mathbb{F}_q^b$. Furthermore, we will identify copies of $\mathcal{T}$ rooted at $\textbf{w}$ with vectors of the form $(x_1,x_2,...x_a)$, where the $x_i$s represent the $a$ non-roots in our copy of $\mathcal{T}$ in the correct order. \\

 When is $(x_1,x_2,...,x_a)$ in $C$?. It is in $C$ iff (1) all the sets of the form $\{x_{j},x_{j+1},...,x_{j+k-1}\}$ and $\{x _{j},x_{j+1},...,x_{j+k-2},w_{l}\}$ that correspond to edges of $\mathcal{T}$ are actually edges in $X$ and (2) $x_i\neq x_j$ whenever those two vertices are in the same part and (3) $x_i\neq w_j$ whenever those two vertices are in the same part. \\

The first condition is equivalent to for all $i,j$, $ f_i(x_{j},x_{j+1},...,x_{j+k-1})=0$ and for all $i,j$,\\ $f_i(x _{j},x_{j+1},...,x_{j+k-2},w_{l})$ $=0$ whenever this corresponds to an edge of $\mathcal{T}$. So the set of $\{x_1,x_2,...,x_a\}$ that satisfy condition (1) form the rational points of a variety $V$, made up of at most $a \cdot b$ equations, each of degree at most $d$.\\

The second and third condition together make up a system of at most $\binom{a}{2}+a \cdot (b-a+k-1)$ complements of linear equations, so the set of $(x_1,x_2,...,x_a)$s that satisfy these conditions is the complement of the rational points of a variety $U$ made up of the product of at most $\binom{a}{2}+a\cdot (b-a+k-1)$ linear equations.\\

We have $C \cong V(\mathbb{F}_q)\backslash U(\mathbb{F}_q)$, where the 2 varieties, $U$ and $V$, both have bounded complexity. We can then split $V$ into a number of irreducible components $V=V_1 \cup V_2 \cup ... \cup V_v$, where $v$ is bounded as a function of the complexity, ie it depends on $a$,$b$,$d$ and $k$ but not on $q$. Then $C \cong V_1(\mathbb{F}_q)\backslash U(\mathbb{F}_q) \cup V_2(\mathbb{F}_q)\backslash U(\mathbb{F}_q) \cup ... \cup V_v(\mathbb{F}_q)\backslash U(\mathbb{F}_q) $. Now for each irreducible component $V_i$, either $V_i\subset U$ (in which case $V_i \backslash U = \emptyset$, so we can ignore this component), or $V_i \cap U$ has dimension strictly smaller than $V_i$. By the Lang-Weil bound ~\cite{LangWeil}, $|V_i(\mathbb{F}_q)|=(1+O(q^{-1/2})) \cdot q^{dim(V_i)}$, while $|V_i(\mathbb{F}_q)\cap U(\mathbb{F}_q)|=(c+O(q^{-1/2})) \cdot q^{dim(V_i\cap U)}$ for some $c$ that is bounded as a function of the complexity (independent of $q$). Therefore when $q$ is large enough, we have $2q^{dim(V_i)} > |V_i(\mathbb{F}_q) \backslash U(\mathbb{F}_q)|\geq \frac{1}{2}q^{dim(V_i)}$. Adding all the pieces up, we have $2v \cdot q^{dim(V)}> |V(\mathbb{F}_q)\backslash U(\mathbb{F}_q)| \geq \frac{1}{2}q^{dim(V)}$. When $dim(V)\geq 1$, this gives us $|V(\mathbb{F}_q)\backslash U(\mathbb{F}_q)| \geq q/2$. Otherwise, $dim(V)= 0$ and $|V(\mathbb{F}_q)\backslash U(\mathbb{F}_q)|< 2v$.\\

Since $v$ is bounded by a function of $a,b,d$ and $k$, we can set $p$ to be twice the upper bound, which ensures that $|V(\mathbb{F}_q)\backslash U(\mathbb{F}_q)|< p$ when $dim(V)= 0$. Now the lemma is proved: we either have $|C|\geq q/2$ or $|C|<p$, as required, (where $p$ is independent of $q$).\\
\begin{flushright}
$\Box$
\end{flushright}

Continuing on with the main proof, we have $\mathbb{P}(|C|\geq p)=\mathbb{P}(|C|\geq q/2)=\mathbb{P}(|C|^s\leq (q/2)^s)$, which by Markov's inequality is $\leq \frac{\mathbb{E}(|C|^s)}{(q/2)^s}$. We now want to calculate $\mathbb{E}(|C|^s)$. Because $\mathcal{T}^{\leq s}$ was defined to be the set of all graphs you can make by taking the union of $s$ copies of $T$ all rooted at the same place, an element of $|C|^s$ corresponds to a copy of an element $H$ in $\mathcal{T}^{\leq s}$ (obtained by taking the union). Also, for every element $H$ in $\mathcal{T}^{\leq s}$, let $\gamma_s(H)$ be the number of ways of expressing it as a union of $s$ copies of $T$. This means that:

\[
 \mathbb{E}(|C|^s)) \leq \sum_{H\in\mathcal{T}^{\leq s}}\gamma_s(H)\cdot \mathbb{E}(|\{ A\in\Hom(H,X) \, : \, r(A)=\underline{w}\}|)\]
 
 To get any further, we will need the following lemma:\\

\begin{mylem}  For any $H\in \mathcal{T}^{\leq s}$, we have $\mathbb{E}(|\{A\in\Hom(H,X) \, : \, r(A)=\underline{w}\}|)) = q^{b \cdot |H|-a \cdot e(H)}$. In other words, the expected number of copies of $H$ rooted at $\underline{w}$ is equal to $q^{b \cdot |H|-a \cdot e(H)}$

\end{mylem}

\paragraph*{Proof:} Call $m=|H|-|R|$. We have: $(x_1,x_2,...,x_m)$ forms a copy of $H$ rooted at $\underline{w}$ if and only if for all $i$, $f_i(x_{j_1},x_{j_2},...,x_{j_k})=0$ whenever this corresponds to an edge of $H$ and for all $i$, $f_i(x_{j_1},x_{j_2},...,x_{j_{k-1},w_{j_k}})=0$ whenever that corresponds to an edge of $H$. The $f_i$s are independent from each other so we only need to find, for each $i$, the probability that $f_i(x_{j_1},x_{j_2},...,x_{j_k})=0$ whenever this corresponds to an edge of $H$ and $ f_i(x_{j_1},x_{j_2},...,x_{j_{k-1},w_{j_k}})=0$ whenever that corresponds to an edge of $H$.\\
\\
For simplicity, we shall call the $e(H)$ points in $(\mathbb{F}_q^b)^k$ corresponding to edges of $H$: $y_1,y_2,...,y_{e(H)}$ and fix them. We want to calculate  $\mathbb{P}(\, \forall j \, \, f_i(y_j)=0)$, knowing that $f_i$ is a random polynomial of degree $\leq d$. We can first without loss of generality make a change of variable $\pi$ such that the first coordinate of each $y_j$ is different. To do so, we proceed as follows: a change of variable is just an invertable $bk\times bk$ matrix acting on the $y_j$s. The first coordinates of $\pi(y_j)$ is given by the dot product of $y_j$ with the first row vector of $\pi$. Given any $j$ and $j'$, the first coordinate of $\pi(y_j)$ is equal to the first coordinate of $\pi(y_{j'})$ if and only if the elements of first row vector of $\pi$ satisfy some linear equation. Thus by repeating this operation over all choices of $j,j'$, we get a set of $\binom{e(H)}{2}$ linear equations in $bk$ variables. The set of all possible first rows for $\pi$ has size $q^{bk}-1$ (we have $bk$ coordinates and the only thing we require is that not all of them are 0). The set of all possible first rows that satisfy one particular linear equation has size $q^{bk-1}-1$ (there is some variable that we can express as a function of the $bk-1$ others, and we still disallow the 0). So if we disallow all first rows that satisfy one of the equations, we end up with at least $q^{bk}-1-\binom{e(H)}{2}(q^{bk-1}-1)$ possible first rows of $\pi$.  Note that because $H\in \mathcal{T}^s$, we have $e(H)\leq sb=d+1$, and since we assumed that $q>\binom{d+1}{2}$, we have $\binom{e(H)}{2}/q<1$. Thus, this number is positive, so there is some choice for a first row of $\pi$ that makes the first coordinate of each $\pi(y_j)$ different. From there, add on the other $bk-1$ rows of $\pi$ arbitrarily just making sure that $\pi$ is invertable. On top of replacing the $y_j$s, we'll also be replacing $f_i$ with $f_i\pi^{-1}$ so that $f_i(y_j)$ stays the same. Note that because $f_i$ was chosen uniformly at random amongst polynomials of degree at most $d$ and because $\pi$ is a bijection, $f_i\pi^{-1}$'s distribution  is also uniform amongst polynomials of degree at most $d$. Therefore without loss of generality, we can assume that the first coordinate of the $y_j$s are distinct. We'll let $z_j$ be the first coordinate of $y_j$.  \\

Now suppose we are given a random polynomial of degree at most $d$ : $f(x_1,x_2,...,x_{kb})$. Consider the coefficients in front of the terms $1$ , $x_1$ , $x_1^2$ , $x_1^3$,... and $x_1^{e(H)}$; call them $c_0$,$c_1$,...,$c_{e(H)}$ respectively. These $c_i$s are random variables chosen independently and uniformly in $\mathbb{F}_q$. We can write $f$ as:

\[f=c_0 + c_1x_1+c_2x_1^2+...+c_{e(H)} x_1^{e(H)}+f'\]

where $f'$ consists of all the other terms that aren't already written down. By letting $c'_{e(H)-1}=c_{e(H)-1}+c_{e(H)}z_{e(H)}$, we can rewrite $c_{e(H)}x_1^{e(H)} + c_{e(H)-1}x_1^{e(H)-1}$ as $c_{e(H)}x_1^{e(H)-1}(x_1-z_{e(H)}) + c'_{e(H)-1}x_1^{e(H)-1}$. Note that since $c_{e(H)-1}$ was chosen uniformly at random in $\mathbb{F}_q$ independent of all the other $c$s and independly of $f'$, $c'_{e(H)-1}$ also has the same properties. We can repeat this process multiple times until we write $f$ as:

\[f = c'_0 + (x_1-z_1)\left[c'_1 + (x_1-z_2)\left[c'_2+ (x_1-z_3)\left[...\left[c'_{e(H)-1} + (x_1-z_{e(H)})c'_{e(H)}\right]...\right]\right]\right]+f'\]

where all the $c'_i$s are uniformly chosen in $\mathbb{F}_q$ independently of each other and independently of $f'$.

Suppose we fix $f'$. The polynomial is 0 at $y_1$ iff $c'_0=-f'(y_1)$, which has probability $1/q$. Then given that  $f(y_1)=0$, the polynomial is 0 at $y_2$ iff $c'_1=\frac{c'_0+f'(y_2)}{z_1-z_2}$, which also has probability $1/q$ (remember that all the $z_i$s were distinct so w are not dividing by 0). We continue in this fashion by induction until we reach $f(y_{e(H)})$ is 0 with probability $1/q$ given that all the others are also 0.Multiplying everything together, we get that the probability that $f(y_j)=0$ for all $j$ is $q^{-e(H)}$.\\
\\
Going back to the last inequality, we get the probabilibity that $(x_1,x_2,...,x_m)$ forms a copy of $H$ rooted at $\underline{w}$ is equal to $\prod_{i=1}^a q^{-e(H)} = q^{-a \cdot e(H)}$. Therefore, the expected number of copies of $H$ rooted at $\underline{w}$ is equal to $q^{b \cdot |H|-a \cdot e(H)}$ and the lemma is proved.\\
\begin{flushright}
$\Box$
\end{flushright}

Remember from Lemma 4 that for all $H$s in $\mathcal{T}^{\leq s}$, we have $e(H)\geq|H|\frac{b}{a}$, so by combining Lemmas 4 and 6 we get: $\mathbb{E}(|\{ A\in\Hom(H,X) \, : \, r(A)=\underline{w}\}|)\leq 1$. 

Putting this back in the previous inequality, we have: 
\begin{eqnarray*}\mathbb{E}(|C|^s) & \leq & \sum_{H\in\mathcal{T}^{\leq s}}\gamma_s(H) \cdot \mathbb{E}(|\{ A\in\Hom(H,X) \, : \, r(A)=\underline{w}\}|) \\
&\leq & \sum_{H\in\mathcal{T}^{\leq s}} \gamma_s(H)
\end{eqnarray*}\\
which is a constant depending only on $s$. We will call this $\beta_s$.\\
\\
Again putting this back into the first inequality, we get: $\mathbb{P}(|C|\geq p)\leq \frac{\mathbb{E}(|C|^s)}{(q/2)^s} \leq$ $\frac{2^s\beta_s}{q^s}$.\\
\\

At this point, we know that when we pick $w_1,w_2,...,w_{b-a+k-1}$ at random (in the correct parts), we have a probability of less than $\frac{2^s\beta_s}{q^s}$ of finding a hypergraph of $\mathcal{T}^p$ rooted at $\textbf{w}$. Let $D$ be the number of choices for $\textbf{w}$ that do lead to finding such a hypergraph. $\mathbb{E}(D) \leq k! \cdot (q^b)^{(b-a+k-1)} \cdot \frac{2^s\beta_s}{q^s}$. But now remember that $s$ was defined as $b(b-a+k-1)+a+1$, so we get $\mathbb{E}(D) \leq \frac{k!2^s\beta_s}{q^{a+1}}$. \\
\\
At this point we're finally ready to reconsider the cases where $X$ is empty. We can split the expectation of $D$ into the case where $X$ is empty and the case where it is not: \\$\mathbb{E}(D)= \mathbb{E}(D \, | \, \text{$X$ empty}).\mathbb{P}(\text{$X$ empty}) + \mathbb{E}(D \, | \, \text{$X$ non-empty}).\mathbb{P}(\text{$X$ non-empty})$.\\ We clearly have no copies of the forbidden hypergraphs when $X$ is empty, so $\mathbb{E}(D \, | \, \text{$X$ empty}) = 0$. Meanwhile, we know from earlier that $\mathbb{P}(\text{$X$ non-empty})\geq q^{-a}$. Putting this together, we get:\\
 $\mathbb{E}(D \, | \, \text{$X$ non-empty}) \leq \mathbb{E}(D).q^a \leq \frac{k!2^s\beta_s}{q}$ \\
\\
Now this has order $\Theta(1/q)$ so when $q$ is large enough, we get $\mathbb{E}(D \, | \, \text{$X$ non-empty}) < 1$ . This proves that there is some choice of $f_1,f_2,...,f_{b-a+k-1}$ for which $X$ is non-empty but that gives no elements of $\mathcal{T}^p$ inside $X$.\\

Thus, we have constructed a hypergraph $X$ with $\Theta(N^{k-\frac{a}{b}})$ edges and that does not contain any element of $\mathcal{T}^p$. Thus, $ex(n,\mathcal{T}^p)= \Omega(n^{k-\frac{a}{b}})$ and the proof of the lower bound is complete.\\
\\

\begin{flushright}
$\Box$
\end{flushright}

\section{The upper bound}

In ~\cite{BukhConlon}, the upper bound used the fact that given graph, we can pick a subgraph with high minimal codegree. This fact is not true in general for hypergraphs, so we will have to do something different. We will instead use a hypergraph version of Sidorenko's conjecture ~\cite{Sidorenko} (also posed by Erdos and Simonovits ~\cite{Simonovits} ) that applies to our hypergraph $\mathcal{T}$. The conjecture states that if $H$ is a $k$-partite $k$-hypergraph with $|H|$ vertices and $e(H)$  edges and if $X$ is a $k$-hypergraph with $n$ vertices and $n^{k-r}/k!$ edges, then there are at least $n^{|H|-e(H)r}$ homomomophic copies of $H$ in $X$. We give the proof from ~\cite{Szegedy} that this is true in the case of hypertrees. However, just having a homomorphism isn't enough; we also need it to be injective, otherwise it could intersect itself. We will use a probabilibistic method to fix this problem and prove that the number of actual (non-self-intersecting) copies of $\mathcal{T}$ in $X$ is also of order at least $\Omega(n^{|T|-e(T)r})$. Finally at the end, we use this large number of copies of $\mathcal{T}$ to find a copy of $\mathcal{T}^p$.\\

\begin{mydef}

In a $k$-hypergraph, a set of vertices is called \emph{adjacent} if there exists a single edge that contains all of them.

\end{mydef}

\begin{mydef} A \emph{hypertree} is the hypergraph version of a tree. They can be constructed, starting from a single edge, by adding edges one at a time such that every new edge intersects the old hypergraph in a set of $k-1$ adjacent vertices. Notice this makes it so there is exactly 1 new vertex for every new edge) \\

\end{mydef}

Also notice that $\mathcal{T}$ is a hypertree.

\paragraph{Entropy method\\}

\begin{mydef}

In a $k$-hypergraph $G$, an \emph{ordered edge} is a sequence of $k$ vertices $(x_1,x_2,...,x_k)$ such that the corresponding $k$-set $\{x_1,x_2,...,x_k\}$ is an edge. Notice how, for each edge, we can find exactly $k!$ corresponding ordered edges.

\end{mydef}

\begin{mylem} 

Let $H$ be a $k$-hypertree with $|H|$ vertices and $e(H)$ edges, and let $X$ be a $k$-hypergraph with $n$ vertices and $n^{k-r}/k!$ edges (so $X$ has $n^{k-r}$ ordered edges). Then $|\Hom(H,X)|\geq n^{|H|-e(H)r}$.

\end{mylem}

\paragraph*{Reformulation:} The statement of the lemma is equivalent to $-\ln(\frac{|\Hom(H,X)|}{n^{|H|}}) \leq -e(H)\ln(n^{-r})$, which is equivalent to $-\ln(\frac{|\Hom(H,X)|}{n^{|H|}}) \leq -e(H)\ln(\frac{|\Hom(E,X|)}{n^{k}})$ ,where $E$ is the hypergraph consisting of just a just a single edge. This is again equivalent to:

\[\ln(n^{|H|}) + \sum_{\Hom(H,X)} \ln\left(\frac{1}{|\Hom(H,X)|}\right)\frac{1}{|\Hom(H,X)|} \leq -e(H)\ln\left(\frac{|\Hom(E,X|)}{n^{k}}\right)\]

The key point here is to notice that $\sum_{\Hom(H,X)} \ln\left(\frac{1}{|\Hom(H,X)|}\right)\frac{1}{|\Hom(H,X)|}$ is the entropy of the uniform probability distribution on the set $\Hom(H,X)$, and that therefore this is minimal amongst all other distributions $\mu$ on $\Hom(H,X)$.\\

Now for any distribution $\mu$ on $\Hom(H,X)$, we set $D(\mu)=\ln(n^{|H|}) + \sum_{A\in \Hom(H,X)} \ln\left(\mu(A)\right)\mu(A)$, and then because the uniform distribution has the smallest entropy, the problem is equivalent to finding some distribution $\mu$ on $\Hom(H,X)$ with 

\[D(\mu)\leq e(H)\ln\left(\frac{n^k}{|\Hom(E,X)|}\right)\]

Two other examples of distributions (which will end up being useful) are:\\

 $\bullet \, \epsilon$, the uniform distribution on ordered edges, where $\epsilon(B)=\frac{1}{n^{k-r}}$ for any ordered edge $B$ of $X$.\\
 
 $\bullet \, \kappa$, the distribution on ordered $(k-1)$ sets where for any ordered $(k-1)$-set $T$ of $X$, we have $\kappa(T)=\frac{deg(T)}{n^{k-r}}$ [$deg(T)$ means the number of edges that contain $T$]. This is a well defined probability distribution because the sum of all the probabilities is equal to 1: $\sum_T \frac{deg(T)}{n^{k-r}} = \sum_T \sum_{B \text{ edge} \, | \, T\subset B} \frac{1}{n^{k-r}} = \sum_{B \text{ edge}} \sum_{T \text{ ordered ($k-1$)-set} \, | \, T \subset B} \frac{1}{n^{k-r}} = \sum_{B \text{ edge}} \frac{k!}{n^{k-r}} = 1$ \\
 
 Using the $\epsilon$ distribution, we can immedidiately reformulate the problem as finding some distribution $\mu$ on $\Hom(H,X)$ such that $D(\mu) \leq e(H)D(\epsilon)$.\\
 
 Some notation: for $A\in \Hom(H,X)$ and $S\subset H$, we write $A_{|S}$ to denote the element of $\Hom(S,X)$ formed by restricting $A$ from $H$ to $S$.

\paragraph{Proof by induction\\}

We will prove by induction on the number of edges, that for any hypertree $H$, there does exist a distribution $\mu$ on $\Hom(H,X)$ such that:\\

(1) $D(\mu)\leq e(H)D(\epsilon)$, and \\

(2) such that for any ordered set $S$ of $k-1$ adjacent vertices of $H$, and for any ordered set $T$ of $k-1$ vertices in $X$, we have: $\sum_{A\in \Hom(H,X): A_{|S}=T} \mu(A) = \kappa(T)=\frac{deg(T)}{n^{k-r}}$.

This means that if you pick an ordered $(k-1)$-set of $X$ at random by first picking a copy of $H$ following distribution $\mu$, and then picking any one set of adjacent $k-1$ vertices , then that will be exactly the same as simply picking an ordered $(k-1)$-set of $X$ at random by following the $\kappa$ rule where each $(k-1)$-set is weighted by its degree. \\
\\

The base case is when $H$ is a single edge, in which case we take $\mu=\epsilon$, which assigns to each edge of $X$ the probability $\frac{1}{n^{k-r}}$. Then (1) holds by definition. As for (2), we have for any $k-1$ set $T$ of $X$, $\sum_{\text{edges that contain $T$}}\frac{1}{n^{k-r}} = \frac{deg(T)}{n^{k-r}}$, as required.\\
\\

Now suppose we have a hypertree $H$ and we want to add another edge $E$ to it to make a bigger hypertree $G$, and where $K=E\cap H$ is the place where the edge is being added. Because we are constructing a hypertree, $K$ is a set of $(k-1)$ adjacent vertices. By the induction hypothesis, we have some distribution $\mu$ on $\Hom(H,X)$ with properties (1) and (2). See the following diagram for reference. \\

 \begin{tikzpicture}[scale=0.85]
 
 \draw (0,0) ellipse (3 and 2) ;
 \draw (2.6,0) circle (1) ;
 \draw (-1,-0.5) circle (0.6) ;
 
 \draw (1.5,2) node {$G$} ;
 \draw (0,1) node {$H$} ;
 \draw (2.4,0) node {$K$} ;
 \draw (3.2,0) node {$E$} ;
 \draw (-1,-0.5) node {$S$} ;
 
 \draw (10,0) ellipse (3 and 2) ;
 \draw (12.6,0) circle (1) ;
 \draw (9,-0.5) circle (0.6) ;
 
 \draw (11.5,2) node {$C$} ;
 \draw (10,1) node {$A$} ;
 \draw (12.4,0) node {$D$} ;
 \draw (13.2,0) node {$B$} ;
 \draw (9,-0.5) node {$T$} ;
 
 \draw (10.2,1) ellipse (3.6 and 4) ;
 \draw (10.2,4) node {$X$} ;
 \draw (5,0) node {\Huge{$\rightarrow$}} ;
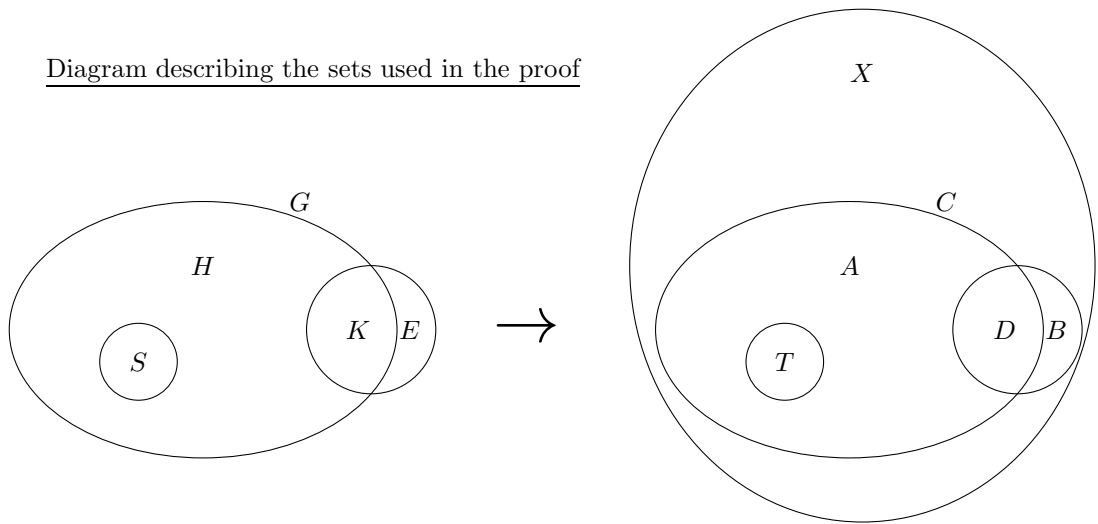
 \draw (1.7,4) node {\underline{Diagram describing the sets used in the proof}} ;

 \end{tikzpicture}

Now we define a probability distribution $\lambda$ on $\Hom(G,X)$ as follows: 

For $C$ a homomorphic copy of $G$, write $C=A\cup B$, where $A$ is the corresponding homomorphic copy of $H$ and $B$ is a single edge corresponding to $E$. Also let $D$ be the set of $(k-1)$ vertices corresponding to $K$ (ie: where $A$ gets attached to $B$). We now define: 

\[\lambda(C)=\mu(A)\frac{1}{deg(D)}\]

We are essentially picking a random copy of $G$ by first picking a random copy of $H$ and then adding another edge to it uniformly at random among all the available choices. Note that there will always be SOME choice (or equivalently, $deg(D)\geq 1$) because $D$ is a set of adjacent vertices, so are contained in an edge.  \\
\\

This distribution $\lambda$ satisfies property (1) because: 

\begin{eqnarray*} & & \ln(n^{|G|})+\sum_{C\in \Hom(G,X)}\ln\left(\lambda(C)\right)\lambda(C) \\
& =  & \ln(n^{|H|+1})+\sum_{\substack{D\in \Hom(K,X) \\ A\in \Hom(H,X) \, : \, A_{|K}=D \\ B\in \Hom(E,X) \, : \, B_{|K}=D}}\ln\left(\mu(A)\frac{1}{deg(D)}\right)\mu(A)\frac{1}{deg(D)} \\
& = & \ln(n^{|H|})+\ln(n)+\sum_{\substack{D \in \Hom(K,X) \\ A\in \Hom(H,X) \, : \, A_{|K}=D}}\ln\left(\mu(A)\frac{1}{deg(D)}\right)\mu(A)\\
& = & \ln(n^{|H|})+\ln(n)+\sum_{\substack{D\in \Hom(K,X) \\ A\in \Hom(H,X) \, : \, A_{|K}=D}}\ln\left(\mu(A)\right)\mu(A)-\sum_{\substack{D\in \Hom(K,X) \\ A\in \Hom(H,X) \, : \, A_{|K}=D}}\ln\left(deg(D)\right)\mu(A)
\end{eqnarray*}

We use property (2) applied to $H,D$ which make this:

\begin{eqnarray*}
& = & \ln(n^{|H|})+\ln(n)+\sum_{A\in \Hom(H,X)}\ln\left(\mu(A)\right)\mu(A)-\sum_{D\in \Hom(K,X)}\ln\left(deg(D)\right)\frac{deg(D)}{n^{k-r}}\\
& = & D(\mu)+\ln(n)-\sum_{D\in \Hom(K,X)}\ln\left(deg(D)\right)\frac{deg(D)}{n^{k-r}}\\
& = & D(\mu) - \left[\ln(n^{k-1}) + \sum_{D\in \Hom(K,X)}\ln\left(\frac{deg(D)}{n^{k-r}}\right)\frac{deg(D)}{n^{k-r}}\right]+\left[\ln(n^k)+\sum_{D\in \Hom(K,X)}\ln\left(\frac{1}{n^{k-r}}\right)\frac{deg(D)}{n^{k-r}}\right]\\
& = & D(\mu)-D(\kappa)+D(\epsilon) \\
& \leq & [e(H)+1]D(\epsilon) \qquad \qquad \text{[$\leftarrow$ we used property (1) applied to $H$]}\\
& = & e(G)D(\epsilon) 
\end{eqnarray*}

Therefore the bigger hypretree $G$ also satisfies property (1) as required.\\
\\
\\
Now we prove $G$ satisfies property (2). Let $S$ be a set of $k-1$ adjacent vertices in $G$. There are two possibilities: either $S$ is a set of adjacent vertices in $H$, or $S$ is contained within the new edge $E$.\\

In the case where $S\subset H$, we have, for any $(k-1)$-set $T$ in $X$, 

\begin{eqnarray*}
& & \sum_{C\in \Hom(G,X) \, : \, C_{|S}=T} \lambda(C) \\
& = & \sum_{\substack{D\in \Hom(K,X) \\ A\in \Hom(H,X) \, : \, A_{|K}=D \,  \text{ and } \,  A_{|S}=T \\ B\in \Hom(E,X) \, : \, B_{|K}=D}}\mu(A)\frac{1}{deg(D)} \\
& = & \sum_{\substack{D\in \Hom(K,X) \\ A\in \Hom(H,X) \, : \, A_{|K}=D \,  \text{ and } \,  A_{|S}=T}}\mu(A) \\
& = & \sum_{A\in \Hom(H,X) \, : \, A_{|S}=T}\mu(A) \\
& = & \frac{deg(T)}{n^{k-r}}  \qquad \qquad \text{[$\leftarrow$ we used property (2) applied to $H,T$]}
\end{eqnarray*}

In the second case where $S\subset E$, we get:

\begin{eqnarray*}
& & \sum_{C\in \Hom(G,X) \, : \, C_{|S}=T} \lambda(C) \\
& = & \sum_{\substack{D\in \Hom(K,X) \\ A\in \Hom(H,X) \, : \, A_{|K}=D \\ B\in \Hom(E,X) \, : \, B_{|K}=D \,  \text{ and } \,  B_{|S}=T}}\mu(A)\frac{1}{deg(D)} \\
& = & \sum_{\substack{D\in \Hom(K,X) \\ B\in \Hom(E,X) \, : \, B_{|K}=D \,  \text{ and } \,  B_{|S}=T}}\frac{1}{n^{k-r}} \qquad \qquad \text{[$\leftarrow$ we used property (2) applied to $H,D$]} \\
& = & \sum_{B\in \Hom(E,X) \, : \, B_{|S}=T}\frac{1}{n^{k-r}} = \frac{deg(T)}{n^{k-r}}
\end{eqnarray*}

So the new distribution $\lambda$ on $\Hom(G,X)$ satisfies property (2).\\
\\
\\
By induction on all hypertrees, we can conclude that for any hypertree $H$, there exists some distribution $\mu$ on $\Hom(H,X)$ such that $D(\mu) \leq e(H)D(\epsilon)$, and therefore all hypertrees satisfy Sidorenko's conjecture; we have $n^{|H|-e(H)r}$ homomorphic copies of $H$ in $X$.\\

\begin{flushright}
$\Box$
\end{flushright}

\paragraph{Probability that a homomorphism is injective\\}

The above is strong indication that we should be able to find a large number of copies of $H$ in $X$. However, in the above we only found homomorphisms from $H$ to $G$. What we actually want is the following:

\begin{mylem}

Let $H$ be a $k$-hypertree with $|H|$ vertices and $e(H)$ edges, and let $X$ be a $k$-hypergraph with $n$ vertices and $n^{k-r}$ edges. Then $|\Inj(H,X)|= \Omega(n^{|H|-e(H)r})$

\end{mylem}

To prove this lemma, we will ask ourselves what is the probability that a homomorphism from $H$ to $X$ (picked at random following distribution $\mu$) is not actually injective? We will bound the answer by using induction.\\
\\
The base case is a single edge $E$, which cannot intersects itself, and therefore the probability is 0.\\
\\
Now suppose $G=H\cup E$ is a hypertree, and $C=A\cup B$ is a copy of our hypertree in $X$, as in the previous part. If $C$ intersects itself, then either $A$ intersects itself or $B$ intersects $H$. ($B$ cannot intersect itself because it's a single edge.) By the induction hypothesis, we know the probability of $A$ intersecting itself, so all that remains is to show that the probability that $B$ intersects with $A$ is not too big.\\

Remember that $\lambda(C)=\mu(A)\frac{1}{deg(D)}$. If we fix $A$ and pick a suitable copy of $E$ at random, then note that adding $E$ only adds 1 single vertex to the hypergraph. Therefore there are at most $|A|-k+1$ possibilities for $E$ to intersect $A$ (each of which corresponds to the new vertex being somewhere in $A$ but outside $A\cap B$). So the probability that $B$ intersects $A$ is at most $\mu(A)\frac{|A|-k+1}{deg(D)}\leq \mu(A)\frac{|H|-k+1}{deg(D)}$. Summing these up across all choices of $A$, we get 

\begin{eqnarray*}
& & \sum_{\substack{D\in \Hom(K,X) \\ A\in \Hom(H,X) \, : \, A_{|K}=D}}\frac{\mu(A)}{deg(D)}(|H|-k+1) \\
& = & \sum_{D\in \Hom(K,X)}\frac{1}{n^{k-r}}(|H|-k+1) \qquad \qquad \text{[$\leftarrow$ we used property (2) applied to $H,D$]} \\
& \leq & \frac{n^{k-1}}{n^{k-r}}(|H|-k+1) \\
& = & n^{r-1}(|H|-k+1)
\end{eqnarray*}

Summing everything up, we get that the overall probability that $G$ intersects itself is at most:\\ $ n^{r-1}[1+2+3+...+(|G|-k)] =  n^{r-1}\frac{|G|^2-(2k-1)|G|+k(k-1)}{2}$. When $r<1$, this is small.\\

For clarity, we will set the probability that a hypertree of size $|G|$ intersects itself to be $\mathbb{P}_{|G|} \leq n^{r-1}\frac{|G|^2-(2k-1)|G|+k(k-1)}{2}$.\\

\paragraph{Size of $\Inj(G,X)$\\}

We know the number of homomorphic images of $G$, and we know the probability that one is injective. From this, we want to find the number of injective images of $G$. This is slightly more complicated than it seems because the probability is not uniform. However, we can still find a lower bound. First, we do another entropy inequality:

\begin{eqnarray*}
\ln(|\Inj(G,X)|) & = & \sum_{C\in \Inj(G,X)} -\ln\left(\frac{1}{|\Inj(G,X)|}\right)\frac{1}{|\Inj(G,X)|} \\
& \geq & \sum_{C\in \Inj(G,X)} -\ln\left(\frac{\lambda(C)}{1-\mathbb{P}_{|G|}}\right)\frac{\lambda(C)}{1-\mathbb{P}_{|G|}} \\
& = & \ln(1-\mathbb{P}_{|G|})+\frac{\sum_{C\in \Inj(G,X)} -\ln(\lambda(C))\lambda(C)}{1-\mathbb{P}_{|G|}}
\end{eqnarray*}

So now we want to know find a lower bound on $\sum_{C\in \Inj(G,X)} -\ln(\lambda(C))\lambda(C)$. We claim that $\sum_{C\in \Inj(G,X)} -\ln(\lambda(C))\lambda(C) \geq \ln\left(n^{(1-r)|G|+(k-1)r}\right) - O\left(n^{r-1}\ln(n)\right)$. We will prove it by induction on $|G|$. For $|G|=k$, ie $G$ is a single edge, we get $\sum_{C\in \Hom(G,X)} -\ln(\lambda(C))\lambda(C) = \sum_{C\in \Hom(G,X)} -\ln\left(\frac{1}{n^{k-r}}\right)\frac{1}{n^{k-r}} = \ln(n^{k-r})= \ln(n^{(1-r)k+(k-1)r})$ so it is true for $|G|=k$. For larger $G$s, we have:

\begin{eqnarray*} 
& & \sum_{C\in \Inj(G,X)} -\ln(\lambda(C))\lambda(C) \\
& = & \sum_{\substack{D\in \Hom(K,X) \\ A\in \Inj(H,X) \, : \, A_{|K}=D \\ B\in \Hom(E,X) \, : \, B_{|K}=D \text{ and } B\cap A=D} } -\ln\left(\frac{\mu(A)}{deg(D)}\right)\frac{\mu(A)}{deg(D)} \\
& \geq & \sum_{\substack{D\in \Hom(K,X) \\ A\in \Inj(H,X) \, : \, A_{|K}=D}} -\ln\left(\frac{\mu(A)}{deg(D)}\right)\mu(A)\frac{deg(D)-|H|}{deg(D)} \\
& = & \left[\sum_{A\in \Inj(H,X)} -ln(\mu(A))\mu(A)\right] +  \left[\sum_{\substack{D\in \Hom(K,X) \\ A\in \Inj(H,X) \, : \, A_{|K}=D}}\ln\left(\frac{deg(D)}{n^{k-r}}\right)\mu(A)\right] \\
& & + \left[\sum_{A\in \Inj(H,X)}\ln(n^{k-r})\mu(A)\right]+\left[\sum_{\substack{D\in \Hom(K,X) \\ A\in \Inj(H,X) \, : \, A_{|K}=D}}\ln\left(\frac{\mu(A)}{deg(D)}\right)\frac{\mu(A)}{deg(D)}|H|\right]
\end{eqnarray*}

There are 4 terms in this inequality, which we will simplify seperately. The first can be found by induction on $|H|$:

\[\sum_{A\in \Inj(H,X)} -\ln(\lambda(A))\lambda(A) \geq \ln\left(n^{(1-r)|H|+(k-1)r}\right) - O\left(n^{r-1}\ln(n)\right) \]

For the second one we have:

\[\sum_{\substack{D\in \Hom(K,X) \\ A\in \Inj(H,X) \, : \, A_{|K}=D}}\ln\left(\frac{deg(D)}{n^{k-r}}\right)\mu(A) \geq \sum_{\substack{D\in \Hom(K,X) \\ A\in \Hom(H,X) \, : \, A_{|K}=D}}\ln\left(\frac{deg(D)}{n^{k-r}}\right)\mu(A)\]

\[= \sum_{D\in \Hom(K,X)}\ln\left(\frac{deg(D)}{n^{k-r}}\right)\frac{deg(D)}{n^{k-r}} \geq \ln\left(\frac{1}{|\Hom(K,X)|}\right) = -\ln(n^{k-1}) \]

For the third term, we have $\sum_{A\in \Inj(H,X)}\ln(n^{k-r})\mu(A)=\ln(n^{k-r})(1-\mathbb{P}_{|H|})$ by definition of $\mathbb{P}_{|H|}$. Now using the fact that $\mathbb{P}_{|H|}\leq n^{r-1}\frac{|H|^2-(2k-1)|H|+k(k-1)}{2}$, we get that this is $\ln(n^{k-r}) - O(\left(n^{r-1}\ln(n)\right)$ Finally, for the fourth term, we have:

\[\sum_{\substack{D\in \Hom(K,X) \\ A\in \Inj(H,X) \, : \, A_{|K}=D}}\ln\left(\frac{\mu(A)}{deg(D)}\right)\frac{\mu(A)}{deg(D)}|H| \geq \sum_{\substack{D\in \Hom(K,X) \\ A\in \Hom(H,X) \, : \, A_{|K}=D}}\ln\left(\frac{\mu(A)}{deg(D)}\right)\frac{\mu(A)}{deg(D)}|H|\]

Remember that we had for all $D\in \Hom(K,X)$,  $\sum_{A\in \Hom(H,X)\, : \, A_{|K}=D}\mu(A) = \frac{deg(D)}{n^{k-r}}$. Using this property, we get that $\sum_{\substack{D\in \Hom(K,X) \\ A\in \Hom(H,X) \, : \, A_{|K}=D}}\frac{\mu(A)}{deg(D)} = \sum_{D\in \Hom(K,X)}\frac{1}{n^{k-r}}=n^{r-1}$. We can now do another entropy inequality to get that the fourth term is at least: 

\begin{eqnarray*}
& & \sum_{\substack{D\in \Hom(K,X) \\ A\in \Hom(H,X) \, : \, A_{|K}=D}}\ln\left(\frac{\mu(A)}{deg(D)}\right)\frac{\mu(A)}{deg(D)}|H| \\
& \geq & \sum_{\substack{D\in \Hom(K,X) \\ A\in \Hom(H,X) \, : \, A_{|K}=D}}\ln\left(\frac{n^{r-1}}{|\Hom(H,X)|}\right)\frac{n^{r-1}}{|\Hom(H,X)|}|H| \\
& = & \ln\left(\frac{n^{r-1}}{|\Hom(H,X)|}\right)n^{r-1}|H|
\end{eqnarray*} 

So now we want to bound $|\Hom(H,X)|$ from above. We can take the trivial bound where we pick an ordered edge, and then $|H|-k$ other vertices arbitrarily. This gives us $|\Hom(H,X)|\leq n^{k-r}n^{|H|-k} = n^{|H|-r}$. Plugging this in gives:

\begin{eqnarray*}
\sum_{\substack{D\in \Hom(K,X) \\ A\in \Hom(H,X) \, : \, A_{|K}=D}}\ln\left(\frac{\mu(A)}{deg(D)}\right)\frac{\mu(A)}{deg(D)}|H| & \geq & \ln\left(\frac{n^{r-1}}{ n^{|H|-r}}\right)n^{r-1}|H| \\
& = & |H|\left(2r-1-|H|)\right)\ln(n)n^{r-1} \\
& = & -O(n^{r-1}\ln(n))
\end{eqnarray*}

Thus, by adding up the four terms back together again, we get: 

\begin{eqnarray*}
\sum_{C\in \Inj(G,X)} -\ln(\lambda(C))\lambda(C) & \geq & \ln\left(n^{(1-r)|H|+(k-1)r}\right) - \ln(n^{k-1}) +\ln(n^{k-r}) - O\left(n^{r-1}\ln(n)\right) \\
& = & \ln\left(\frac{n^{(1-r)|H|+(k-1)r}.n^{k-r}}{n^{k-1}}\right) - O\left(n^{r-1}\ln(n)\right) \\
& = & \ln\left(n^{(1-r)|G|+(k-1)r}\right) - O\left(n^{r-1}\ln(n)\right)
\end{eqnarray*}

This completes the induction. Remember from earlier that \\ $\ln(|\Inj(G,X)|)\geq \ln(1-\mathbb{P}_{|G|})+\frac{\sum_{C\in \Inj(G,X)} -\ln(\lambda(C))\lambda(C)}{1-\mathbb{P}_{|G|}}$. Replacing both\\ $\sum_{C\in \Inj(G,X)} -\ln(\lambda(C))\lambda(C)\geq \ln\left(n^{(1-r)|G|+(k-1)r}\right) - O\left(n^{r-1}\ln(n)\right)$ and \\ $\mathbb{P}_{|G|}\leq O\left(n^{r-1}\right)$, we get:

\begin{eqnarray*}
& & \ln(|\Inj(G,X)|) \\
& \geq & \ln\left(1- O\left(n^{r-1}\right)\right) + \frac{\ln\left(n^{(1-r)|G|+(k-1)r}\right) - O\left(n^{r-1}\ln(n)\right)}{1- O\left(n^{r-1}\right)} \\
& = & \ln\left(1- O\left(n^{r-1}\right)\right) + \ln\left(n^{(1-r)|G|+(k-1)r}\right) - O\left(n^{1-r}\ln(n)\right)+\ln\left(n^{(1-r)|G|+(k-1)r}\right)*O\left(n^{r-1}\right) \\
& = & \ln\left(1- O\left(n^{r-1}\right)\right) + \ln\left(n^{(1-r)|G|+(k-1)r}\right) - O\left(n^{r-1}\ln(n)\right)+O\left(n^{r-1}\ln(n)\right) \\
& = & \ln\left(1- O\left(n^{r-1}\right)\right) + \ln\left(n^{(1-r)|G|+(k-1)r}\right) + \ln\left(1- O\left(n^{r-1}\ln(n)\right)\right) \\
& = & \ln\left(\left(n^{(1-r)|G|+(k-1)r}\right)\left(1-O\left(n^{r-1}\right)\right)\left(1-O\left(n^{r-1}\ln(n)\right)\right)\right) \\
& = & \ln\left(\left(n^{(1-r)|G|+(k-1)r}\right)\left(1-O\left(n^{r-1}\ln(n)\right)\right)\right)
\end{eqnarray*}

Thus, $|\Inj(G,X)|\geq \left(n^{(1-r)|G|+(k-1)r}\right)\left(1-O\left(\frac{\ln(n)}{n^{1-r}}\right)\right)$.

\begin{flushright}
$\Box$
\end{flushright}

\paragraph{Finishing up the upper bound\\}

Our hypergraph $\mathcal{T}$ from the first part is a hypertree, so we know that there are at least $\left(n^{(1-r)|\mathcal{T}|+(k-1)r}\right)\left(1-c\frac{\ln(n)}{n^{1-r}}\right) = \left( n^{b+k-1-rb}\right)\left(1-c\frac{\ln(n)}{n^{1-r}}\right)$ copies of it in any large enough hypergraph $X$ with $n$ vertices and $n^{k-r}$ ordered edges (where $c$ is a constant depending only on $|\mathcal{T}|$). Now note that there are only $n^{b-a+k-1}$ possibilities for choosing distinct roots. Therefore, given a random ordered set of $b-a+k-1$ vertices of $X$, the expected number of copies of $\mathcal{T}$ rooted at them is $\frac{\left( n^{b+k-1-rb}\right)\left(1-c\frac{\ln(n)}{n^{1-r}}\right)}{n^{b-a+k-1}} = n^{a-rb}\left(1-c\frac{\ln(n)}{n^{1-r}}\right)$.\\

If $r$ is small (ie: there are lots of edges), then this quantity will be bigger than any constant. An element of $\mathcal{T}^{\leq p-1}$ can have at most $(p-1)a$ non-roots, which means that we can find at most $[(p-1)a]!/[(p-2)a]!$ copies of $\mathcal{T}$ in it (just choose the order of the vertices). So if our $r$ is small enough that we have at least $[(p-1)a]!/[(p-2)a]!+1$ copies of $\mathcal{T}$ all rooted at the same place, then their union cannot fit into $\mathcal{T}^{\leq p-1}$, therefore it has to be an element of $\mathcal{T}^{\leq u}\backslash \mathcal{T}^{\leq p-1}$ for some large $u$. By removing extra edges and vertices as needed, we can find an element of  $\mathcal{T}^{\leq p}\backslash \mathcal{T}^{\leq p-1} = \mathcal{T}^{p}$. For simplicity, set $p'=[(p-1)a]!/[(p-2)a]!+1$.\\
\\
If we set $r=\frac{a}{b}-\frac{\ln(2p')}{b\ln(n)}$, we get that the expected number of copies of $\mathcal{T}$ rooted at a specific place is at least $2p'\left(1-c\frac{\ln(n)}{n^{1-r}}\right)$. Now as $n\rightarrow \infty$, $\frac{\ln(n)}{n^{1-r}}\rightarrow 0$, so in particular, when $n$ is large enough, we will have $\left(1-c\frac{\ln(n)}{n^{1-r}}\right) \geq 1/2$, and hence the expected number of copies of $\mathcal{T}$ rooted at a specific place is at least $p'$. This is an expectation, so there must exist some choice of $b-a+k-1$ roots such that there are actually at least $p'$ copies of $\mathcal{T}$ rooted at them, and hence an element of $\mathcal{T}^{p}$\\

Therefore, if there are at least $(2p')^{1/b}n^{k-\frac{a}{b}}$ edges in the hypergraph $X$ and $n$ is sufficiently large, then we have a copy of a hypergraph of $\mathcal{T}^p$ inside $X$. Therefore $ex(n,\mathcal{T}^p) \leq O(n^{k-\frac{a}{b}})$. Combining this with the lower bound we proved in the first section we get $ex(n,\mathcal{T}^p)= \Theta(n^{k-\frac{a}{b}})$ and the proof is complete.\\
\\
\begin{flushright}
$\Box$
\end{flushright}

\section{The case where $r\geq 1$}

We will now try to prove Theorem 2, which is the generalisation of Theorem 1 from $0\leq r <1$ to $0\leq r <k-1$. To do so, we use the following lemma:

\begin{mylem}

Given a set of $l$-hypergraphs $\mathcal{F}$ (each of which contains 2 disjoint edges) and some $k>l$, there exists some set $\mathcal{F}'$ of $k$-hypergraphs with \\$ex(n+k-l,\mathcal{F}')= ex(n,\mathcal{F})$ for all $n$.

\end{mylem}

 This will directly prove Theorem 2 when we set $l=\lceil k-r \rceil$, because in that case, $k-r=l-r'$ for some $0\leq r'<1$ so we know from Theorem 1 that we have a set of $l$-hypergraphs $\mathcal{F}$ with $ex(n,\mathcal{F})=\Theta(n^{l-r'})=\Theta(n^{k-r})$. Also for $b\geq k$ (which we can assume without loss of generality), there will be at least 2 disjoint edges in every hypergraph. Applying the lemma gives us some set $\mathcal{F}'$ of $k$-hypergraphs with $ex(n,\mathcal{F}')=\Theta(n^{k-r})$. \\

\paragraph*{Proof of Lemma 9:} For a $l$-hypergraph $F$ and vertices $x_1,x_2,...,x_{k-l}$, define the $k$-hypergraph $(F,x_1,x_2,...,x_{k-l})$ to have vertices $F\cup\{x_1,x_2,...,x_{k-l}\}$ and edges $\{E\cup\{x_1,x_2,...,x_{k-l}\}: E \text{ an edge of $F$}\}$. We define $(\mathcal{F},x_1,x_2,...,x_{k-l})=\{(F,x_1,x_2,...,x_{k-l}) : F\in \mathcal{F}\}$\\$ \cup\{\text{$k$-hypergraphs with $\leq l+2$ edges that are not of the form $(H,x_1,x_2,...,x_{k-l})$}\}$. We claim that $\mathcal{F}'=(\mathcal{F},x_1,x_2,...,x_{k-l}))$ will solve the problem.\\

Suppose $G$ is a $l$-hypergraph with $ex(n,\mathcal{F})$ edges that doesn't contain any element of $\mathcal{F}$. Consider $(G,y_1,y_2,...,y_{k-l})$ for some vertices $y_1,y_2,...,y_{k-l}$. Then first of all, every set of $\leq l+2$ edges of $(G,y_1,y_2,...,y_{k-l})$ is of the form $(H,y_1,y_2,...,y_{k-l})$ because every edge contains $y_1,y_2,...,y_{k-l}$. Now suppose it contains a different element of $(\mathcal{F},x_1,x_2,...,x_{k-l})$, say $(F,x_1,x_2,...,x_{k-l})$. Now $F$ contains 2 edges that do not intersect, so $(F,x_1,x_2,...,x_{k-l})$ contains two edges that intersect only at $x_1,x_2,...,x_{k-l}$. Since any two edges of $(G,y_1,y_2,...,y_{k-l})$ intersect at $y_1,y_2,...,y_{k-l}$, that must mean that $y_1,y_2,...,y_{k-l}$ are the images of $x_1,x_2,...,x_{k-l}$ in some order. Now taking that copy of $(F,x_1,x_2,...,x_{k-l})$ in $(G,y_1,y_2,...,y_{k-l})$ and removing $y_1,y_2,...,y_{k-l}$ from all the edges, we end up with a copy of $F$ inside $G$. This contradicts our original assumption about $G$. Therefore $(G,y_1,y_2,...,y_{k-l})$ is a $k$-hypergraph with $ex(n,\mathcal{F})$ edges that does not contain any element of $(\mathcal{F},x)$, so $ex(n+k-l,(\mathcal{F},x_1,x_2,...,x_{k-l}))\geq ex(n,\mathcal{F})$, as required.\\

Respectively, suppose that $G$ is a $k-hypergraph$ that doesn't contain any element of $(\mathcal{F},x_1,x_2,...,x_{k-l})$. Since every set of $l+2$ edges is of the form $(H,x_1,x_2,...,x_{k-l})$, that means that the entire hypergraph is of the form $(K,x_1,x_2,...,x_{k-l})$ for some $K$.\\
 Indeed, suppose for a contradiction that the hypergraph is not of that form. Pick two edges $e_1$ and $e_2$. They intersect in at most $k-1$ places. Pick $k-l$ of those and call the set $S$. Now the hypergraph is not of the form $(H,x_1,x_2,...,x_{k-l})$, so that means that there is some other edge, $e_3$, that doesn't contain $S$. Now $e_1$,$e_2$,$e_3$ intersect in at most $k-2$ places. Repeat the argument sevral times until we get edges $e_1$, $e_2$, $e_3$, ..., $e_{l+2}$ that intersect in at most $k-l-1$ places. Thus, we have $l+2$ edges that are not of the form $(H,x_1,x_2,...,x_{k-l})$, contradicting our assumption. Therefore $G$ must be of the form $(H,x_1,x_2,...,x_{k-l})$. \\
 Then this $K$ cannot contain any element $F$ of $\mathcal{F}$ because otherwise $G=(K,x_1,x_2,...,x_{k-l})$ would contain $(F,x_1,x_2,...,x_{k-l})$. Therefore $K$ has at most $ex(n,\mathcal{F})$ edges, so $G$ also has at most $ex(n,\mathcal{F})$ edges. Thus $ex(n+k-l,(\mathcal{F},x_1,x_2,...,x_{k-l}))\leq ex(n,\mathcal{F})$\\

 This proves that $ex(n,\mathcal{F})=ex(n+k-l,(\mathcal{F},x_1,x_2,...,x_{k-l}))$, completing the proof of the lemma, and thus proving Theorem 2.\\
 \begin{flushright}
$\Box$
\end{flushright}

\paragraph*{Remarks:\\}

\subparagraph*{The case $k>r>k-1$ is impossible}:\\ Suppose that $\mathcal{F}$ is a collection of $k$-graphs which has $ex(n,\mathcal{F})=\Theta(n^{k-r})$ for some $k>r>k-1$. \\

Now consider $X$ to be the $k$-hypergraph with $n$ vertices defined as follows: it consists of some set $S$ of $t$ vertices, for some $0\leq t\leq k-1$. The other $n-t$ vertices are partitioned into $\lfloor(n-t)/(k-t) \rfloor$ sets of size $k-t$, which we will call $e_1,e_2,...,e_{\lfloor(n-t)/(k-t) \rfloor}$. The edges of the hypergraph are exactly $e_i\cup S$ for $1\leq i\leq \lfloor(n-t)/(k-t) \rfloor$. This hypergraph $X$ has the property that the intersection of any two edges is exactly $S$ therefore it is a sunflower. It also has $\Theta(n)$ edges, which is larger than $ c \cdot n^{k-r}$ for large enough $n$. Therefore $\mathcal{F}$ must contain a subgraph of $X$. However, any subgraph of $X$ must also have the property that any the intersection of two edges is exactly $S$, ie it is another sunflower. We will call this sunflower $F_t$\\

In this way, we get for all $0\leq t \leq k-1$, a sunflower $F_t$ in $\mathcal{F}$ with kernal size $t$. The sunflower lemma states that when this occurs, $ex(n,\mathcal{F})$ has order $O(1)$. This contradicts our assumtion that $ex(n,\mathcal{F})= \Theta(n^{k-r})$. Therefore it is indeed impossible to have a collection of $k$-hypergraphs with $ex(n,\mathcal{F})=\Theta(n^{k-r})$ for any $k>r>k-1$.

\subparagraph{The case $r=k$ is possible:} The hypergraph $E$ consisting of just 1 edge has $ex(n,E)=0$. \\

\subparagraph{The case $r=k-1$ is possible:} Firstly, when $k=2$, it is fairly easy to see that the path of length 3, $P_3$, has $ex(n,P_3)=n$ or $n-1$, which is of the correct order $\Theta(n)$. For larger $k$, we simply apply Lemma 9 to get a collection of $k$-hypergraphs $\mathcal{F}$ with $ex(n,\mathcal{F})=\Theta(n)$ as required. \\
\\

So in conclusion, the rationals $r$ for which there exist some finite $\mathcal{F}$ with $ex(n,\mathcal{F})=\Theta(n^{k-r})$ are exactly  those in the set: $\{r\in \mathbb{Q} : 0\leq r \leq k-1\}\cup\{k\}$.

\paragraph{Acknowledgements:\\}

Thanks to my supervisor, Oleg Pikhurko, for help and guidance.

\bibliography{ex_n_F}{}
\bibliographystyle{plain}
\end{document}